\documentclass[12pt]{article}
\usepackage{amssymb,amsmath,amsfonts,enumerate}
\title{
\author{\bf{Lech Pasicki}}
\bf{Nonexpansive fixed point theorems for primitive uniform spaces}
\footnote{MSC:  54H25}}
\newtheorem{theorem}{\indent Theorem}
\newtheorem{lemma}[theorem]{\indent Lemma}
\newtheorem{proposition}[theorem]{\indent Proposition}
\newtheorem{definition}[theorem]{\indent Definition}
\newtheorem{corollary}[theorem]{\indent Corollary}

\newcommand{\mN}{\mbox{$\mathbb{N}$}}
\newcommand{\FXX}{\mbox{$F\colon X \to 2^{X}$}}
\newcommand{\FmY}{\mbox{$F_{\mid Y}$}}
\newcommand{\fmY}{\mbox{$f_{\mid Y}$}}
\newcommand{\FnY}{\mbox{$F^{n}(Y)$}}
\newcommand{\fnY}{\mbox{$f^{n}(Y)$}}
\newcommand{\fnpY}{\mbox{$f^{n+1}(Y)$}}
\newcommand{\cB}{\mbox{$c_{\cal{B}}$}}
\newcommand{\eyx}{\mbox{$ \emptyset \neq Y \subset X$}}
\newcommand{\Bpi}{\mbox{$\cal{B}$}}
\newcommand{\XB}{\mbox{$(X,\cal{B})$}}
\newcommand{\Upi}{\mbox{$\cal{U}$}}
\newcommand{\ai}{\mbox{$\alpha_{i}$}}
\newcommand{\si}{\mbox{$x_{i}$}}
\newcommand{\yi}{\mbox{$y_{i}$}}
\newcommand{\co}{\mbox{$conv \;$}}
\begin{document}
\maketitle
\vspace{1 in}
\begin{abstract}
Some known fixed point theorems for nonexpansive mappings in metric spaces are extended here  
to the case of primitive uniform spaces. The reasoning presented in the proofs seems to be a natural way to 
obtain other general results.
\end{abstract}

\par In \cite{pa1}, \cite{pa2}, \cite{pa3} some fixed point theorems for
metric spaces were proved.  It is a natural idea to consider a more general case. 
\begin{definition}
\label{De1}
   Let $X$ be a nonempty set, and let $\Bpi$ be a family of subsets of $X \times X$
such that $X \times X \in \Bpi$. Then $\Bpi$ is a primitive uniformity for $X$ and
$\XB$ is a primitive uniform space.
\end{definition}
\begin{definition}[{\rm cp. \cite[Definition 5]{pa1}\bf}]
\label{De2}
Let $\XB$ be a primitive uniform space, and let $A$ be a nonempty subset of $X$.  An $x \in X$ is a 
central point for $A$ if
\begin{equation}
\label{co1}
  \bigcap \{B \in \Bpi: A \subset B(z) \mbox{ for a } z \in X \}=
  \bigcap \{B \in \Bpi: A \subset B(x) \}
\end{equation}
holds;  the centre for $A$ is the set $c_{\Bpi}(A)$ of all such central points.
\end{definition}
\par In particular, if for a metric space $(X,d)$ we adopt the natural base 
\begin{equation*}
  \Bpi=\{\{(x,y) \in X \times X\colon d(x,y) <r\}\colon r > 0 \} \cup \{X \times X\},
\end{equation*}
then $c_{\Bpi}(A)$  coincides with $c(A)$ defined in \cite[Definition 5]{pa1}.
\par Let $2^X$ be the family of all subsets of $X$, and then $\FXX$ is a (multivalues) mapping 
if $F(x) \neq \emptyset$, $x \in X \neq \emptyset$.
\begin{definition}[{\rm cp. \cite[Definition 7]{pa1}\bf}]
\label{De3}
Let $\XB$ be a primitive uniform space, $\eyx$, and $F\colon Y \to 2^Y$ a mapping. An $x \in X$ is a central 
point for $F$ if
\begin{equation}
\label{co2}
\begin{split}
  & \bigcap \{B \in \Bpi: \FnY \subset B(z) \mbox{ for an } n \in \mN
          \mbox{ and a } z \in Z \}=\\
  & \bigcap \{B \in \Bpi: \FnY \subset B(x)
	    \mbox{ for an } n \in \mN \}
\end{split}
\end{equation}
holds; the centre $c_{\Bpi}(F)$ for $F$ is the set of all such central points.
\end{definition}
\par  In particular, for a metric space $(X,d)$, $\eyx$, and $F\colon Y \to 2^{Y}$,
we adopt the following (see \cite[(8)]{pa3})
\begin{equation}
\label{co3}
\begin{split}
  & r(F):=inf\{t \in (0,\infty]: \FnY \subset B(z,t)
     \mbox{ for a } z \in X \mbox{ and} \\   
  & \mbox{an } n \in \mN \} = 
      inf \{t \in (0,\infty]: \FnY \subset B(x,t) \mbox{ for an } n \in \mN \} \mbox{.}
\end{split}
\end{equation}
The set $ c(F) $ consists then of points $x$ satisfying \eqref{co3}.
\par The subsequent extension of \cite[Theorem 13]{pa3} transfers the problem of existence of fixed points to 
the problem of finding spaces (and bases) for which centres are singletons.
\begin{theorem}
\label{Th4}
Let $\XB$ be a primitive uniform space, and let $f \colon X \to X$ be a mapping.
Assume that $\eyx$ is such that $f(Y) \subset Y$, $x \in \cB(\fmY)$ and the following
\begin{equation}
\label{co4}
    f(Y \cap B(x)) \subset B(f(x)), \quad B \in \Bpi
\end{equation}
holds. Then $f(x) \in \cB(\fmY)$, and $x$ is a fixed point for $f$ provided that $\cB(\fmY)$ is a singleton.
\end{theorem}
\par Proof.  We have $\fnpY \subset \fnY \subset \cdots \subset f(Y) \subset Y$.
If $\fnY \subset B(x)$ holds, then by \eqref{co4} we obtain 
\[
\fnpY = f(\fnY \cap B(x)) \subset f(Y \cap B(x)) \subset B(f(x)),
\]
i.e. $f(x) \in \cB(\fmY)$. $\quad \square$
\par In \cite[Lemma 13]{pa2} it was stated that if $(X,d)$ is a bead space
(\cite[Definition 6]{pa2}, then for any mapping $\FXX$ and nonempty bounded set $Y \subset X$ 
such that $F(Y) \subset Y$, the set $\cB(\FmY)$ is at most a singleton. Bead spaces are extensions of
uniformly convex spaces and of unitary spaces.
\par The subsequent definition shows that condition \eqref{co4} is quite a natural one.
\begin{definition}
\label{De5}
Let $\XB$ be a primitive uniform space. Then a mapping $F\colon X \to 2^X$ is nonexpansive, if the following 
is satisfied
\begin{equation}
\label{co5}
   F(B(x)) \subset B(F(x)), \quad x \in X \mbox{, } B \in \Bpi \mbox{.}
\end{equation}
\end{definition}
\par Clearly, for $F\colon X \to X$ condition \eqref{co5} is a particular
case of condition \eqref{co4} and therefore the theorem to follow is a consequence of
Theorem \ref{Th4}.
\begin{theorem}
\label{Th6}
Let $\XB$ be a primitive uniform space, and let $f \colon X \to X$ be a nonexpansive mapping. Assume that 
$\eyx$ is such that $f(Y) \subset Y$. Then $x \in \cB(\fmY)$ yields $f(x) \in \cB(\fmY)$, and 
and $x$ is a fixed point for $f$ provided that $\cB(\fmY)$ is a singleton.
\end{theorem}
\par It is clear that $\cB(\fmY)$ need not be singletons.  If we are interested in fixed point theorems 
for mappings satisfying condition \eqref{co4} or \eqref{co5}, then we should investigate spaces 
$\XB$ and subsets $Y$ of $X$. 
\par In what follows, if $\{x\}=\cB(F(z))$, then we adopt $(\cB \circ F)(z)=x$.
\par The next theorem is a consequence of Theorem \ref{Th4}, and it extends \cite[Theorem 15]{pa3}.
\begin{theorem}
\label{Th7}
Let $\XB$ be a primitive uniform space, and let $\FXX$ be a mapping such that $\cB(F(x))$ are singletons with 
$\cB(F(x)) \subset F(x)$, $x \in X$.  Assume that $\eyx$ is such that for $f = \cB \circ F$ we have 
$f(Y) \subset Y$ (e.g. if $F(Y) \subset Y$) and $x \in \cB(\fmY)$. If condition \eqref{co4} 
is satisfied, then $f(x) \in \cB(\fmY)$, and $x$ is a fixed point for $F$ provided that $\cB(\fmY)$ is a singleton.
\end{theorem}
\par Proof. In view of Theorem \ref{Th4} we have $f(x) \in \cB(\fmY)$, and 
$\{x\} = \{f(x)\} = \cB(F(x)) \subset F(x)$ provided that $\cB(\fmY)$ is a singleton. $\quad \square$
\par The following is a consequence of Theorems \ref{Th6}, \ref{Th7}.
\begin{theorem}
\label{Th8}
Let $\XB$ be a primitive uniform space, and let $\FXX$ be a mapping such that $\cB(F(x))$ are singletons with 
$\cB(F(x)) \subset F(x)$, $x \in X$. Assume that $f=\cB \circ F$ is nonexpansive, and for a $\eyx$ 
we have $f(Y) \subset Y$ (e.g. if $F(Y) \subset Y$) and $x \in \cB(\fmY)$. Then $f(x) \in \cB(\fmY)$, and $x$ is 
a fixed point for $F$ provided that $\cB(\fmY)$ is a singleton.
\end{theorem}
\par  Assume that $(X,d)$ is a metric space.
Condition 
\begin{equation*}
   F(B(x,r)) \subset B(F(x),r), \quad x \in X 
\end{equation*}
(cp. \eqref{co5}) means that
\begin{equation*}
   F(y) \subset B(F(x),r), \quad y \in B(x,r), 
\end{equation*}
and we obtain
\begin{corollary}
\label{Cor9}
Let $(X,d)$ be a metric space, and let $\FXX$ be a mapping. Then
$F$ is nonexpansive in the sense of Definition \ref{De5} (for the 
natural base in X) if and only if
\begin{equation}
  \label{co6}
  D(F(x),F(y)) \leq d(x,y), \quad x,y \in X
\end{equation}
holds, where $D$ is the Hausdorff metric.
\end{corollary}
\par It is much easier to prove the following (certainly known) proposition, than to find it :).
\begin{proposition}
\label{Pr10}
Assume that $Z$ is a linear space, $A,C,V \subset Z$ are nonempy and $V$ convex. Then
$A \subset V+C$ yields $\co A \subset V+\co C$.
\end{proposition}
\par Proof.  Let $ x_{1},\ldots,x_{n} \in A$ be arbitrary, and let
$ y_{1},\ldots,y_{n} \in C $ be such that $\si \in V+ \yi$. Then for $ x= \sum \ai \si$
(a convex combination) we obtain $x \in \sum \ai (V+ \yi) = V+(\sum \ai \yi)$ and, consequently, 
$\co A \subset V+ \co C$. $\quad \square$
\par If $X$ is a convex set in a linear topological space $Z$, then it can be treated as 
a uniform space $(X,\Upi)$, where for any $U \in \Upi$ and all $x \in X$
we have $U(x) = (V+x) \cap X$, for the respective neighbourhood $V$ of zero in $Z$.
\begin{lemma}
\label{Le11}
Let Z be a locally convex space, and let $(X,\Upi)$ be its uniform  subspace with $X$ convex. Assume that 
$\FXX$ is a mapping, for a $\eyx$ we have $F(Y) \subset Y$, and \eqref{co4} holds with $F$ in place of 
$f$ for a base $\Bpi$ of $\Upi$ with convex $B(s)$, $B \in \Bpi$, $s \in X$. Then the following is satisfied
\begin{equation}
\label{co7}
   (\co \circ F)(Y \cap B(x)) \subset B(\co F(x)), \quad B \in \Bpi \mbox{.}
\end{equation}
\end{lemma}
\par Proof. We have
\begin{equation*}
   (\co \circ F)(Y \cap B(x)) = \bigcup \{ \co F(z): z \in Y \cap B(x) \} \subset 
   \co F(Y \cap B(x)).
\end{equation*}
For each $B \in \Bpi$ there exists a convex set $V$ such that $B(C) = V + C$. Now, we apply Proposition 
\ref{Pr10} to $A=F(Y \cap B(x))$, and $C = F(x)$ while $A=F(Y \cap B(x)) \subset B(F(x)) = V + F(x)$ 
(see \eqref{co4}). $\quad \square$
\begin{corollary}
\label{Cor12}
Let Z be a locally convex space, and let $(X,\Upi)$ be its 
uniform  subspace with $X$ convex. Assume that $\FXX$ is a nonexpansive mapping
with respect to a base $\Bpi$ of $\Upi$ with convex $ B(x)$, $B \in \Bpi$,
$x \in X$. Then $\co \circ F$ is nonexpansive with respect to $\Bpi$.
\end{corollary}
\begin{lemma}
\label{Le13}
Let $\XB$ be a primitive uniform space, and $\eyx$. If $f\colon X \to X$ is a mapping such that 
$f(Y) \subset Y$, and \eqref{co4} holds for all $x \in Z = c_{\cal{B}}(\fmY)$, then  
$f_{\mid Z}\colon Z \to Z$. Assume that $F\colon Y \to 2^Y$ is a mapping. If $(X,d)$ is a metric space, 
then $c(F)$ is closed; if $X$ is a convex set in a linear space, all $B \in \Bpi$ are symmetric, and 
$B(x) \cap X$ is convex  for each $B \in \Bpi$, $x \in X$, then $c_{\Bpi}(F) \cap X$ is convex.
\end{lemma}
\par Proof. If $x \in Z$ then in view of Theorem \ref{Th4} $f(x) \in Z$. Let $(X,d)$ be a metric space, and let 
us adopt $Z=c(F)$, $b = r(F) < \infty$ (see \eqref{co3}). If  $x \in \overline{Z}$, then for any $\epsilon > 0$
there exists a $z \in Z$ such that $d(x,z) < \epsilon$. Consequently, we obtain 
$\FnY \subset B(z,b + \epsilon) \subset B(x,b + 2\epsilon)$ for large $n$, which means that $x \in Z$.  If 
$x,z \in \cB(F)$ and $B \in \Bpi$ is symmetric, then $y \in \FnY \subset B(x) \cap B(z)$ yields 
$x,z \in B(y)$. Moreover, if $x,z \in \cB(F) \cap X$ and $B(y) \cap X$ is convex, then we have 
$\alpha x + (1-\alpha)z \in B(y) \cap X$ and $y \in B(\alpha x + (1-\alpha)z)$ for any $\alpha \in [0,1]$ ($B$ is 
symmetric), which means that $\alpha x + (1-\alpha)z \in \cB(F)$, and $\cB(F) \cap X$ is convex. $\quad \square$

\section*{Acknowledgements}
This work was partially supported by the Faculty of Applied Mathematics AGH UST statutory tasks within subsidy 
of the Polish Ministry of Science and Higher Education, grant no. 16.16.420.054.
\par 
\vspace{.1in}
\mbox{Faculty of Applied Mathematics} \linebreak
\mbox{AGH University of Science and Technology} \linebreak
\mbox{Al. Mickiewicza 30} \linebreak
\mbox{30-059 KRAK\'OW, POLAND} \linebreak
\mbox{E-mail: pasicki@agh.edu.pl}

\end{document}